\ifx\shlhetal\undefinedcontrolsequence\let\shlhetal\relax\fi
%
\input amstex
\NoBlackBoxes
\loadmsbm 
\def\lesdot{\mathrel{\mathord{<}\!\!\raise 0.8 pt\hbox{$\scriptstyle\circ$}}}  

\def\hang{\hangindent\parindent}
\def\item{\par\hang\textindent}
\def\itemitem{\par\indent \hangindent2\parindent \textindent}

\def\textindent#1{\indent\llap{#1\enspace}\ignorespaces}
\def\narrower{\advance\leftskip by\parindent \advance\rightskip
 by\parindent}
\documentstyle {amsppt}
\topmatter
\title {Possibly every real function is continuous on a non-meagre set} 
\endtitle
\author {Saharon Shelah \thanks {\null\newline
Partially supported by the Edmund Landau Center for Research in 
Mathematical Analysis (supported by the Minerva Foundation, Germany) and the 
``The Israel Science Foundation'' administered by the Israeli Academy
of Sciences and Humanities. Publication 473.\null\newline 
We would like to thank Alice Leonhardt for the excellent typing.\null\newline 
} \endthanks} \endauthor
\affil {Institute of Mathematics \\
The Hebrew University \\
Jerusalem, Israel
\medskip
Rutgers University \\
Mathematics Department \\
New Brunswick, NJ  USA} \endaffil
\endtopmatter
\document
\expandafter\ifx\csname bib4plain.tex\endcsname\relax
  \expandafter\gdef\csname bib4plain.tex\endcsname{}
\else \message{Hey!  Apparently you were trying to \string twice.   This does not make sense.}
\errmessage{Please edit your file (probably \jobname.tex) and remove
any duplicate ``\string\input'' lines} \fi

\def\renewcommand{\newcommand}	       
\edef\cite{\the\catcode`@}%
\catcode`@ = 11
\let\@oldatcatcode = \cite
\chardef\@letter = 11
\chardef\@other = 12
%
%
%
%
\def\@innerdef#1#2{\edef#1{\expandafter\noexpand\csname #2\endcsname}}%
%
%
\@innerdef\@innernewcount{newcount}%
\@innerdef\@innernewdimen{newdimen}%
\@innerdef\@innernewif{newif}%
\@innerdef\@innernewwrite{newwrite}%
%
%
%
\def\@gobble#1{}%
%
%
%
\ifx\inputlineno\@undefined
   \let\@linenumber = \empty 
\else
   \def\@linenumber{\the\inputlineno:\space}%
\fi
%
%
%
\def\@futurenonspacelet#1{\def\cs{#1}%
   \afterassignment\@stepone\let\@nexttoken=
}%
\begingroup 
\def\\{\global\let\@stoken= }%
\\ 
\endgroup
\def\@stepone{\expandafter\futurelet\cs\@steptwo}%
\def\@steptwo{\expandafter\ifx\cs\@stoken\let\@@next=\@stepthree
   \else\let\@@next=\@nexttoken\fi \@@next}%
\def\@stepthree{\afterassignment\@stepone\let\@@next= }%
%
%
%
\def\@getoptionalarg#1{%
   \let\@optionaltemp = #1%
   \let\@optionalnext = \relax
   \@futurenonspacelet\@optionalnext\@bracketcheck
}%
%
%
\def\@bracketcheck{%
   \ifx [\@optionalnext
      \expandafter\@@getoptionalarg
   \else
      \let\@optionalarg = \empty
      \expandafter\@optionaltemp
   \fi
}%
\def\@@getoptionalarg[#1]{%
   \def\@optionalarg{#1}%
   \@optionaltemp
}%
%
%
%
\def\@nnil{\@nil}%
\def\@fornoop#1\@@#2#3{}%
\def\@for#1:=#2\do#3{%
   \edef\@fortmp{#2}%
   \ifx\@fortmp\empty \else
      \expandafter\@forloop#2,\@nil,\@nil\@@#1{#3}%
   \fi
}%
\def\@forloop#1,#2,#3\@@#4#5{\def#4{#1}\ifx #4\@nnil \else
       #5\def#4{#2}\ifx #4\@nnil \else#5\@iforloop #3\@@#4{#5}\fi\fi
}%
\def\@iforloop#1,#2\@@#3#4{\def#3{#1}\ifx #3\@nnil
       \let\@nextwhile=\@fornoop \else
      #4\relax\let\@nextwhile=\@iforloop\fi\@nextwhile#2\@@#3{#4}%
}%
%
%
%
\@innernewif\if@fileexists
\def\@testfileexistence{\@getoptionalarg\@finishtestfileexistence}%
\def\@finishtestfileexistence#1{%
   \begingroup
      \def\extension{#1}%
      \immediate\openin0 =
         \ifx\@optionalarg\empty\jobname\else\@optionalarg\fi
         \ifx\extension\empty \else .#1\fi
         \space
      \ifeof 0
         \global\@fileexistsfalse
      \else
         \global\@fileexiststrue
      \fi
      \immediate\closein0
   \endgroup
}%
%
%
%
%
\def\bibliographystyle#1{%
   \@readauxfile
   \@writeaux{\string\bibstyle{#1}}%
}%
\let\bibstyle = \@gobble
%
%
\let\bblfilebasename = \jobname
\def\bibliography#1{%
   \@readauxfile
   \@writeaux{\string\bibdata{#1}}%
   \@testfileexistence[\bblfilebasename]{bbl}%
   \if@fileexists
      \nobreak
      \@readbblfile
   \fi
}%
\let\bibdata = \@gobble
%
%
\def\nocite#1{%
   \@readauxfile
   \@writeaux{\string\citation{#1}}%
}%
\@innernewif\if@notfirstcitation
%
%
\def\cite{\@getoptionalarg\@cite}%
%
%
\def\@cite#1{%
   \let\@citenotetext = \@optionalarg
   \printcitestart
   \nocite{#1}%
   \@notfirstcitationfalse
   \@for \@citation :=#1\do
   {%
      \expandafter\@onecitation\@citation\@@
   }%
   \ifx\empty\@citenotetext\else
      \printcitenote{\@citenotetext}%
   \fi
   \printcitefinish
}%
\def\@onecitation#1\@@{%
   \if@notfirstcitation
      \printbetweencitations
   \fi
   \expandafter \ifx \csname\@citelabel{#1}\endcsname \relax
      \if@citewarning
         \message{\@linenumber Undefined citation `#1'.}%
      \fi
      \expandafter\gdef\csname\@citelabel{#1}\endcsname{%
\strut
\vadjust{\vskip-\dp\strutbox
\vbox to 0pt{\vss\parindent0cm \leftskip=\hsize 
\advance\leftskip3mm
\advance\hsize 4cm\strut\openup-4pt 
\rightskip 0cm plus 1cm minus 0.5cm ?  #1 ?\strut}}
         {\tt
            \escapechar = -1
            \nobreak\hskip0pt
            \expandafter\string\csname#1\endcsname
            \nobreak\hskip0pt
         }%
      }%
   \fi
   \csname\@citelabel{#1}\endcsname
   \@notfirstcitationtrue
}%
%
%
\def\@citelabel#1{b@#1}%
%
%
\def\@citedef#1#2{\expandafter\gdef\csname\@citelabel{#1}\endcsname{#2}}%
%
%
%
\def\@readbblfile{%
   \ifx\@itemnum\@undefined
      \@innernewcount\@itemnum
   \fi
   \begingroup
      \def\begin##1##2{%
         \setbox0 = \hbox{\biblabelcontents{##2}}%
         \biblabelwidth = \wd0
      }%
      \def\end##1{}
      %
      %
      \@itemnum = 0
      \def\bibitem{\@getoptionalarg\@bibitem}%
      \def\@bibitem{%
         \ifx\@optionalarg\empty
            \expandafter\@numberedbibitem
         \else
            \expandafter\@alphabibitem
         \fi
      }%
      \def\@alphabibitem##1{%
         \expandafter \xdef\csname\@citelabel{##1}\endcsname {\@optionalarg}%
         \ifx\biblabelprecontents\@undefined
            \let\biblabelprecontents = \relax
         \fi
         \ifx\biblabelpostcontents\@undefined
            \let\biblabelpostcontents = \hss
         \fi
         \@finishbibitem{##1}%
      }%
      \def\@numberedbibitem##1{%
         \advance\@itemnum by 1
         \expandafter \xdef\csname\@citelabel{##1}\endcsname{\number\@itemnum}%
         \ifx\biblabelprecontents\@undefined
            \let\biblabelprecontents = \hss
         \fi
         \ifx\biblabelpostcontents\@undefined
            \let\biblabelpostcontents = \relax
         \fi
         \@finishbibitem{##1}%
      }%
      \def\@finishbibitem##1{%
         \biblabelprint{\csname\@citelabel{##1}\endcsname}%
         \@writeaux{\string\@citedef{##1}{\csname\@citelabel{##1}\endcsname}}%
         \ignorespaces
      }%
      %
      %
      \let\em = \bblem
      \let\newblock = \bblnewblock
      \let\sc = \bblsc
      \frenchspacing
      \clubpenalty = 4000 \widowpenalty = 4000
      \tolerance = 10000 \hfuzz = .5pt
      \everypar = {\hangindent = \biblabelwidth
                      \advance\hangindent by \biblabelextraspace}%
      \bblrm
      \parskip = 1.5ex plus .5ex minus .5ex
      \biblabelextraspace = .5em
      \bblhook
      \input \bblfilebasename.bbl
   \endgroup
}%
%
%
\@innernewdimen\biblabelwidth
\@innernewdimen\biblabelextraspace
%
%
%
\def\biblabelprint#1{%
   \noindent
   \hbox to \biblabelwidth{%
      \biblabelprecontents
      \biblabelcontents{#1}%
      \biblabelpostcontents
   }%
   \kern\biblabelextraspace
}%
%
%
%
\def\biblabelcontents#1{{\bblrm [#1]}}%
%
%
\def\bblrm{\rm}%
%
%
\def\bblem{\it}%
%
%
\def\bblsc{\ifx\@scfont\@undefined
              \font\@scfont = cmcsc10
           \fi
           \@scfont
}%
%
%
\def\bblnewblock{\hskip .11em plus .33em minus .07em }%
%
%
\let\bblhook = \empty
%
%
%
\def\printcitestart{[}
\def\printcitefinish{]}
\def\printbetweencitations{, }
\def\printcitenote#1{, #1}
%
%
%
\let\citation = \@gobble
%
%
%
\@innernewcount\@numparams
%
%
\def\newcommand#1{%
   \def\@commandname{#1}%
   \@getoptionalarg\@continuenewcommand
}%
%
%
\def\@continuenewcommand{%
   \@numparams = \ifx\@optionalarg\empty 0\else\@optionalarg \fi \relax
   \@newcommand
}%
%
%
\def\@newcommand#1{%
   \def\@startdef{\expandafter\edef\@commandname}%
   \ifnum\@numparams=0
      \let\@paramdef = \empty
   \else
      \ifnum\@numparams>9
         \errmessage{\the\@numparams\space is too many parameters}%
      \else
         \ifnum\@numparams<0
            \errmessage{\the\@numparams\space is too few parameters}%
         \else
            \edef\@paramdef{%
               \ifcase\@numparams
                  \empty  No arguments.
               \or ####1%
               \or ####1####2%
               \or ####1####2####3%
               \or ####1####2####3####4%
               \or ####1####2####3####4####5%
               \or ####1####2####3####4####5####6%
               \or ####1####2####3####4####5####6####7%
               \or ####1####2####3####4####5####6####7####8%
               \or ####1####2####3####4####5####6####7####8####9%
               \fi
            }%
         \fi
      \fi
   \fi
   \expandafter\@startdef\@paramdef{#1}%
}%
%
%
%
%
\def\@readauxfile{%
   \if@auxfiledone \else 
      \global\@auxfiledonetrue
      \@testfileexistence{aux}%
      \if@fileexists
         \begingroup
            \endlinechar = -1
            \catcode`@ = 11
            \input \jobname.aux
         \endgroup
      \else
         \message{\@undefinedmessage}%
         \global\@citewarningfalse
      \fi
      \immediate\openout\@auxfile = \jobname.aux
   \fi
}%
%
%
\newif\if@auxfiledone
\ifx\noauxfile\@undefined \else \@auxfiledonetrue\fi
%
%
%
%
\@innernewwrite\@auxfile
\def\@writeaux#1{\ifx\noauxfile\@undefined \write\@auxfile{#1}\fi}%
%
%
%
\ifx\@undefinedmessage\@undefined
   \def\@undefinedmessage{No .aux file; I won't give you warnings about
                          undefined citations.}%
\fi
%
%
\@innernewif\if@citewarning
\ifx\noauxfile\@undefined \@citewarningtrue\fi
%
%
%
\catcode`@ = \@oldatcatcode

\newpage

\head {\S0 Introduction} \endhead
\bigskip

By Abraham, Rubin and Shelah \cite{ARSh:153} it is consistent that every
function from ${\Bbb R}$ to ${\Bbb R}$ is continuous when restricted to some
uncountable set (and more ...).  We may consider strengthening this statement,
by demanding the subset on which the function is continuous to be large in a
stronger sense. In \cite[Problem AR(b)]{Fe94}, David Fremlin asked exactly
this, namely, is it consistent that every function from ${\Bbb R}$ to ${\Bbb
R}$ is continuous when restricted to some set which is non-meagre (i.e. not
countable union of nowhere dense sets). We answer it (positively) here. We use
${}^\omega 2$ for reals, and for $B\subseteq {}^\omega 2$ we say $f:B
\rightarrow {}^\omega 2$ is continuous if $f(\eta_0) = \eta_1, n < \omega$
implies that for some $m$ we have  
$$(\forall \eta)[\eta_0\restriction m \triangleleft \eta \in B\Rightarrow
\eta_1 \restriction n \triangleleft f(\eta)],$$
where for sequences $\nu,\eta$, $\nu\triangleleft\eta$ means ``$\eta$
(properly) extends $\nu$''.

The non-meagre sets here are of cardinality $\aleph_1$ and we may wonder
whether we can get them of higher cardinality.  From the point of view of
those who asked the original question, maybe this does not add much (but I
think it does make for much more interesting partition relations).  Another
interesting generalization, asked by Heinrich von Weizs\"acker (see
\cite[Problem AR(a)]{Fe94}), is:
\medskip
\roster
\item "{$(*)$}"   is it consistent to have that every function from
${\Bbb R}$  to  ${\Bbb R}$  is continuous when restricted to
some non-null set? 
\endroster
\medskip

\noindent
We may ask about $2$-place functions; Sierpinski colouring implies we cannot
ask for one colour, but we may want the consistency of:
\medskip
\roster
\item "{$(**)$}"  for every $2$-place function $f:{\Bbb R}\times{\Bbb R}
\longrightarrow {\Bbb R}$ there is a non-meagre set (non-null) $A \subseteq
{\Bbb R}$ and $2$-place continuous functions $f_0,f_1:A\times A\longrightarrow
{\Bbb R}$ such that 
$$(\forall x,y \in A)(f(x,y) = f_0(x,y) \vee f(x,y) = f_1(x,y))$$
\endroster
\medskip

\noindent
(similarly for larger $n$ or all $n$ simultaneously). It is not clear to us
how much this interests non-logicians, but to us it seems to be the right
question and we shall deal with $(**)$ and the cardinality in Rabus Shelah
\cite{RaSh:585}. 

We can also generalize the proof replacing $\aleph_0$ by $\mu = \mu^{<\mu}$
(as done here).  On consistent partition relations see \cite{Sh:276},
\cite{Sh:288}, \cite{Sh:546} and \cite{RaSh:585}.  Another related theorem,
proved in \cite{Sh:481}, is the consistency of
\medskip
\roster
\item "{$(***)$}"  if $B_1$ is a $2^{\aleph_0}$-c.c. Boolean algebra and $B_2$
is a c.c.c. Boolean algebra, then $B_1*B_2$ is a $2^{\aleph_0}$-c.c.
Boolean algebra.
\endroster
\bigskip

Originally the proof goes through ``meagre preserving" and iterations as in
section 2 of \cite{Sh:276}, but the proof was cumbersome, had flaws and proper
version was not manufactured.  Here we present a simpler though a
``degenerated'' variant (i.e. we use only the forcing we have to and so we get
$2^{\aleph_0} = \aleph_2$ rather than $2^{\aleph_0} = \aleph_3$). 
\medskip

We thank Andrzej Ros{\l}anowski and Zoran Spasojevi\'c for much help in
proofreading. 
\bigskip

\head {\S1 Continuity on a Non-Meagre Set} \endhead
\bigskip

\definition{1.1 Definition}  1) Cohen$_\mu(\alpha) = \{f:f$ a partial function
from $\alpha$ to $\{ 0,1\}$ \newline

$\qquad \qquad \qquad \qquad \qquad \qquad \qquad \quad$ 
of cardinality $< \mu\}$ \newline
ordered by inclusion.\newline
2) Cohen$_\mu = \text{ Cohen}_\mu(\mu)$. \newline
3) A $\mu$-Cohen forcing $P$ means Cohen$_\mu(\alpha)$ for some ordinal
$\alpha$ (or at least the set
$$\{ p \in P:P \restriction \{q:q \ge p\}\text{\ is\ equivalent\ to\ some\
Cohen}_\mu(\alpha)\}$$
is a dense subset of $P$).
\enddefinition
\bigskip

\definition{1.2 Definition}  A set $A \subseteq {}^\mu 2$ is $\mu$-meagre if
it is the union of $\le \mu$ nowhere dense subsets of ${}^\mu 2$.
\enddefinition
\bigskip

\proclaim{1.3 Main Lemma}  Let $\mu$ be a regular cardinal such that
$\mu=\mu^{<\mu}$. For $\ell < 2$ and $\alpha<\mu^+$, let $Q_{\alpha,\ell}$
be $({}^{\mu >} 2,\triangleleft)$ (i.e. Cohen$_\mu$ forcing) and for $I
\subseteq\mu^+ \times 2$ let $P_I$ be the product $\prod\limits_{t \in I} Q_t$
with $(<\mu)$-support and $P = P_{\mu^+ \times 2}$. Let ${\underset\tilde {}
\to\eta_{\alpha,\ell}}$ be the $Q_{\alpha,\ell}$-name of the generic function
from $\mu$ to $2 = \{0,1\}$. In 
$$V^P = V[ \langle {\underset\tilde {}\to \eta_{\alpha,\ell}}:
\alpha < \mu^+,\ell < 2 \rangle]$$ 
let $R$ be a forcing notion defined by:
$$
\align
R = \biggl\{ (u_0,v,\bar \nu):&(a) \quad u_0 \text{ is a subset of } \mu^+
\text{ of cardinality } <\mu,\\ 
  &(b) \quad v \text{ is a subset of }{}^{\mu>}2\text{ of cardinality }<\mu\\
  &(c) \quad \text{if } \rho \in v\text{ then } \ell g(\rho) \text{ is a
successor ordinal and}\\
  &\qquad \quad \rho(\ell g(\rho)-1) = 1,\\
  &(d) \quad \bar \nu = \langle \nu_\rho:\rho \in v\rangle \text{ with }
\nu_\rho \in {}^{\mu >} 2,\\
  &(e)  \quad \text{if } \rho_1 \triangleleft \rho_2 \text{ are from } v 
\text{ then } \nu_{\rho_1} \triangleleft \nu_{\rho_2} \text{ (hence they are
not equal)},\\ 
  &(f) \quad \text{for every } \rho \in v\text{ there is } \alpha \in u_0
\text{ such that } \rho\triangleleft {\underset\tilde {}\to\eta_{\alpha,0}},\\ 
  &(g) \quad \text{if } \alpha \in u_0,\rho \triangleleft {\underset\tilde
{}\to\eta_{\alpha,0}}\text{ and } \rho \in v \text{ then } \nu_\rho
\triangleleft {\underset\tilde {}\to\eta_{\alpha,1}}.\biggr\}.
\endalign
$$
\medskip

\noindent
\underbar{The order is natural}:\newline
$(u_0,v,\bar\nu) \le (u'_0,v',\bar \nu')$ if and only if $u_0\subseteq u'_0$,\ 
$v\subseteq v'$, and $\bar\nu=\bar\nu'\restriction v$ (and both are in $R$).

Let ${\underset\tilde {}\to {\Cal U}}$ be an $R$-name such that
$$
\Vdash_R {\underset\tilde {}\to {\Cal U}} = \bigcup\{u_0:\text{there are } 
v,\bar \nu\text{ such that }(u_0,v,\bar\nu)\in{\underset\tilde {}\to G_R}\},
$$
where ${\underset\tilde {}\to G_R}$ is the canonical $R$-name for the generic
filter on $R$. \newline
\underbar{Then}:
\medskip
\roster
\item "{$(a)$}"  In $V^{P * \underset\tilde {}\to R}$, the function
${\underset\tilde{}\to\eta_{\alpha,0}}\mapsto{\underset\tilde{}\to
\eta_{\alpha,1}}$ from $\{{\underset\tilde{}\to\eta_{\alpha,0}}:\alpha\in
{\underset\tilde{}\to{\Cal U}}\}$ to $\{{\underset\tilde{}\to\eta_{\alpha,1}}:
\alpha\in{\underset\tilde{}\to{\Cal U}}\}$ is continuous.
\item "{$(b)$}"   In $V^{P * \underset\tilde {}\to R}$, the set
$\{ {\underset\tilde {}\to \eta_{\alpha,0}}:\alpha \in
{\underset\tilde {}\to {\Cal U}} \}$ is a non-$\mu$-meagre subset of 
${}^\mu 2$, in fact everywhere non-meagre (i.e. its intersection with any
non-empty open set is non-meagre).
\item "{$(c)$}"  $P * \underset\tilde {}\to R$ is a $\mu$-Cohen forcing,
moreover $P_{\alpha \times 2} \lesdot P$ (i.e. $P_{\alpha\times 2}$ is a
complete suborder of $P$) and $P * \underset\tilde{}\to R/P_{\alpha \times 2}$
is a $\mu$-Cohen forcing for every $\alpha < \mu^+$, $\alpha\geq \mu$.
\endroster
\endproclaim
\bigskip

\remark{1.3A Remark}  We can get slightly more than continuity for $\mu>
\aleph_0$ (a slight change of the forcing notion), namely, for some club $W
\subseteq \mu$, for every $i\in W$ we have that ${\underset\tilde {}\to
\eta_{\alpha,0}} \restriction i$ determines ${\underset\tilde {}\to
\eta_{\alpha,1}} \restriction i$ (and so can get this in the theorem). 
\endremark
\bigskip

\demo{Proof}  For $\alpha \le \mu^+,\alpha \ge \mu$ let
$$
\align
A'_\alpha =: \biggl\{(p,r):&(\alpha) \quad p \in P_{\alpha \times 2}
\text{ and} \\
  &\qquad \quad p = \langle \eta^p_{\gamma,\ell}:
(\gamma,\ell) \in u^p \times 2 \rangle \text{ where } \eta_{\gamma,\ell}^p 
\in {}^{\mu >}2, u^p\in [\alpha]^{<\mu},\\
  &(\beta) \quad r = (u_0,v,\bar \nu) \text{ satisfies clauses }
(a), (b), (c), (d), (e) \\
  &\qquad \quad \text{of the definition of } R \text{ and} \\
  &\qquad \quad(f)' 
\quad \text{ for every } \rho \in v \text{ there is } \beta \in
u^p \text{ such that} \\
  &\qquad \qquad \qquad \beta \in u_0 \text{ and } \rho \trianglelefteq
\eta^p_{\beta,0}, \\
  &\qquad \quad(g)'\quad\text{ if } \beta \in u_0 \text{ and } \rho \in v
\text{ and } \rho \trianglelefteq \eta^p_{\beta,0}, \\
  &\qquad\qquad\qquad\text{then }\nu_\rho\trianglelefteq\eta^p_{\beta,1},\\
  &(\gamma)\quad u_0\subseteq u^p.\biggr\}
\endalign
$$
\medskip

\noindent and partially order $A'_\alpha$ in the natural way.

Also let
$$
\align
A_\alpha = \biggl\{ (p,r):&p,r\text{ satisfy clauses }(\alpha),(\beta),
(\gamma) \text{ above and}\\ 
  &(\delta) \quad \text{ if } \beta\in u_0\text{ and } \rho\in v\text{ then }
\neg(\eta^p_{\beta,0} \triangleleft \rho),\\ 
  &(\varepsilon) \quad \text{ for every } \beta \in u_0\text{ there is }\gamma
\in u_0\cap \mu \text{ such that}\\ 
  &\qquad \quad (\eta^p_{\beta,0},\eta^p_{\beta,1}) = (\eta^p_{\gamma,0},
\eta^p_{\gamma,1}) \biggr.\}.
\endalign
$$
\medskip

Note that in $A^\prime_\alpha$ any increasing sequence of length $<\mu$ has a
least upper bound, i.e.~if a sequence $\langle q_i: i<\delta\rangle\subseteq
A^\prime_\alpha$ is increasing, $\delta<\mu$ then there is a condition $q\in
A^\prime_\alpha$ such that $q_i\leq q$ for all $i<\delta$ and if
$\bigwedge\limits_{i<\delta} q_i\leq q'\in A^\prime_\alpha$ then $q\leq q'$. 
\medskip

\noindent Now:
\medskip
\roster
\item "{$(A)$}"  $A_\alpha$ is a dense subset of $A'_\alpha$. If $(p,r) \in
A_\alpha$ and $\mu\leq\alpha \le \mu^+$ then $(p,r) \in P * \underset\tilde
{}\to R$ (i.e. $p \Vdash_P ``r \in  \underset\tilde {}\to R"$).
\endroster
[Why?  The second phrase should be clear. Trivially $A_\alpha \subseteq
A'_\alpha$. Suppose $(p,r) \in A'_\alpha$ and we shall find $(p',r) \in
A_\alpha$ such that $(p,r) \le (p',r)$.  Let $\beta^* =: \sup(u^p \cap \mu) +
1$ and $\varepsilon = \sup\{ \ell g(\rho) + 1:\rho \in v\}$ and let
${\bar{\bold 0}}_\varepsilon$ be the sequence of zeroes of length
$\varepsilon$ and let $u^p=\{\gamma_\zeta:\zeta < |u^p|\}$ be an enumeration.  

Let $u^{p'} = u^p \cup \{\beta^* + \zeta:\zeta < |u^p|\}$ and
$$
\eta^{p'}_{\gamma,\ell} = \cases \eta^p_{\gamma,\ell} \char 94 
{\bar{\bold 0}}_\varepsilon \quad &\text{ \underbar{if} } \quad \gamma \in
u^p,\ell=0, \\ 
  \eta^p_{\gamma,\ell} \quad &\text{ \underbar{if} } \quad 
\gamma \in u^p,\ell=1, \\
  \eta^p_{\gamma_\zeta,\ell} \char 94 {\bar{\bold 0}}_\varepsilon \quad 
&\text{ \underbar{if} } \quad \gamma=\beta^* + \zeta, \ell=0, \\
  \eta^p_{\gamma_\zeta,\ell}\quad &\text{ \underbar{if} }\quad \gamma=\beta^*
+ \zeta,\ell = 1. \endcases 
$$
\medskip
\noindent Note: if $\rho \in v$ then $\rho \trianglelefteq \eta \char 94
{\bar{\bold 0}}_\varepsilon \Leftrightarrow \rho \trianglelefteq \eta$
(because by clause (c) in the definition of $R$, the last element in $\rho$ is
one).] 
\roster
\item "{$(B)$}"  $A_{\mu^+}$ is a dense subset of 
$P * \underset\tilde {}\to R$.
\endroster
[Why?  Since $P$ is $\mu$-closed, for each $(p,\underset\tilde {}\to r) \in P
* \underset\tilde {}\to R$ we can find $p_1 \in P$ and $r$ such that $p \le
p_1 \in P$ and $p_1 \Vdash_P ``\underset\tilde {}\to r = r"$. Next we choose
$p_2$ such that $p_1 \le p_2 \in P$ and clauses $(\alpha)$--$(\varepsilon)$ in
the definition of $A_{\mu^+}$ hold.  So $(p,\underset\tilde {}\to r) \le
(p_2,r)$ (in $P * \underset\tilde {}\to R)$ and $(p_2,r) \in A_{\mu^+}$.]
\medskip

\noindent 
If $\mu \le \alpha \le \beta \le \mu^+$ and $(p,r) \in A_\beta$ then we let
$(p,r) \restriction \alpha = (p \restriction \alpha,r \restriction \alpha)$
where $p \restriction \alpha = \langle \eta^p_{\gamma,\ell}:(\gamma,\ell)\in
(u^p \cap \alpha) \times 2 \rangle$ and $r \restriction \alpha = (u^r_0 \cap
\alpha,v^r,\bar \nu^r)$. It is easy to check that:
\medskip
\roster
\item "{$(C)$}"  If $\mu \le \alpha \le \beta \le \mu^+$ and $(p,r) \in
A_\beta$ then $(p,r) \restriction \alpha \in A_\beta$ and $(p,r) \restriction
\alpha \in A_\alpha$ 
\endroster
(note that clause $(f)'$ of the definition of $A_\alpha$ follows from clause
$(\varepsilon)$ in the definition of $A_\alpha$). 
\roster
\item "{$(D)$}"  If $\mu \le \alpha \le \beta \le \mu^+$ and $(p_1,r_1) \in
A_\beta$ and $(p_1,r_1) \restriction \alpha \le (p_2,r_2) \in A_\alpha$
then $(p_1,r_1),(p_2,r_2)$ are compatible (in $A_\beta$); moreover, they have
the natural least upper bound
$$
\biggl( p_2 \cup p^*,(u^{r_1}_0 \cup u^{r_2}_0,v^{r_2},\bar \nu^{r_2}) \biggr)
$$
in $A'_\beta$, where $p^*\in P_{\beta\times 2}$ is such that
$u^{p^*}=u^{p_1}\setminus\alpha$,
and for $\gamma\in u^{p^*}$
$$
\eta^{p^*}_{\gamma,\ell}=\cases \eta^{p_1}_{\gamma,\ell} &\text{ if }\qquad
\text{ either }\gamma\notin u^{r_1}_0\text{ or }\ell=0,\\
\eta^{p_1}_{\gamma,1}\cup\bigcup\{\nu^{r_2}_\rho: \rho\trianglelefteq
\eta^{p_1}_{\gamma,0}\} &\text{ if }\qquad\gamma\in u^{r_1}_0\text{ and }
\ell=1.
\endcases
$$
\endroster
Moreover:
\roster
\item "{$(D)^+$}" {$\underline{if}$} $\langle (p_{1,i},r_{1,i}):i\leq
\delta\rangle$ is an increasing continuous sequence of conditions in
$A^\prime_\beta$ such that $\delta<\mu$ and $(p_{1,i},r_{1,i})\in A_\beta$ for
each non-limit $i\leq\delta$ and $(p_{1,i},r_{1,i})\restriction\alpha\leq
(p_2,r_2)$ for $i<\delta$ 

\noindent {$\underline{then}$} $(p_{1,\delta},r_{1,\delta})$, $(p_2,r_2)$ are
compatible in $A^\prime_\beta$ (and $(p_{1,\delta},r_{1,\delta})\restriction
\alpha\leq (p_2,r_2)$. 
\endroster
Hence
\roster
\item "{$(D)^{++}$}" if $\mu\leq\alpha\leq\beta\leq\mu^+$ then $A_\alpha
\subseteq A_\beta$, moreover $A_\alpha\lesdot A_\beta$ (by $(D)$) and the
quotient $A_\beta/A_\alpha$ is $\mu$-closed (i.e.~increasing chains of length
$<\mu$ have upper bound). 
\endroster
\remark{Remark} For better understanding note that if $Q_0$ is adding a
$\mu$-Souslin tree, $\underset\tilde {}\to Q_1$ is the forcing determined by
this tree, then the composition $Q_0*\underset\tilde {}\to Q_1$ is not even
$\mu$-closed. The point is that increasing sequences in $Q_0*\underset\tilde
{}\to Q_1$ of length $<\mu$ have an upper bound, but not necessarily least
upper bound. 
\endremark

\noindent Now (note the support)
\roster
\item "{$(E)$}"  the sequence $\langle A_\alpha:\alpha \in [\mu,\mu^+]
\rangle$ is increasing, and for limit $\alpha\in (\mu,\mu^+)$, $A_\alpha$ is
the inverse limit of $\langle A_\beta:\beta < \alpha \rangle$ if cf$(\alpha)<
\mu$, and direct limit of it if cf$(\alpha) \ge \mu$. Hence $\mu\le\alpha<
\beta \le \mu^+ \Rightarrow A_\beta/A_\alpha$ is a $\mu$-Cohen forcing
(remember that atomless $\mu$-closed forcing notion of size $\mu$ is
$\mu$-Cohen). 
\endroster
It is also clear that
\roster
\item "{$(F)$}"  $P_{\alpha \times 2} \lesdot P_{\beta \times 2}$ if
$\alpha < \beta \le \mu^+$.
\endroster
Lastly note that
\roster
\item "{$(G)$}"  for $\alpha \in [\mu,\mu^+),A_\alpha/P_{\alpha \times 2}$
is $\mu$-closed and hence it is $\mu$-Cohen.
\endroster
\medskip

\relax From the desirable conclusions, we have gotten clause $(c)$ (by
combining clauses $(E)$ and $(G)$). 
\medskip
\roster
\item "{$(H)$}"  If $(p_1,r_1)\in A_{\mu^+}$, $\beta^*\in\mu^+\setminus
u^{p_1}$ and $\nu^*\in {}^{\mu>} 2$ then there is $(p_2,r_2)\in A_{\mu^+}$
such that
$$
(p_1,r_1)\leq (p_2,r_2),\quad \beta^*\in u^{r_2}_0\quad\text{ and }\quad
\nu^*\triangleleft \eta^{p_2}_{\beta^*,0}.
$$
\endroster
[Why?  Let $\beta^* \in \mu^+\setminus u^{p_1}$, $\beta^{**}\in \mu\setminus
(u^{p_1}\cup\{\beta^*\})$ and let 
$$
\nu^\otimes = \bigcup\{ \nu_\rho^{r_1}:\rho \trianglelefteq \nu^* \text{ and }
\rho\in v^{r_1}\}\quad\text{ and }\quad\gamma = \sup\{ \ell g(\nu)+1:\nu\in
v^{r_1}\cup\{ \nu^*\}\}.
$$
Now define $(p_2,r_2)$ as follows:
$$
u^{p_2} = u^{p_1} \cup \{ \beta^*,\beta^{**}\}
$$

$$
\eta^{p_2}_{\beta,\ell} =\cases \eta^{p_1}_{\beta,\ell} &\text{ if }
\qquad \beta \in u^{p_1},\ell < 2, \\
  \nu^* \char 94 {\bar{\bold 0}}_\gamma &\text{ if } \qquad \beta \in 
\{ \beta^*,\beta^{**}\},\ell = 0,\\
  \nu^\otimes &\text{ if } \qquad \beta \in
\{\beta^*,\beta^{**}\},\ell=1, 
 \endcases
$$
and $r_2 = (u^{r_2}_0,v^{r_2},\bar \nu^{r_2})$ is defined by $u^{r_2}_0=
u^{r_1}_0\cup \{ \beta^*,\beta^{**}\}$,\quad $v^{r_2}=v^{r_1}$,\quad
$\bar{\nu}^{r_1}=\bar{\nu}^{r_2}$. Easily $(p_2,r_2)$ is as required.]
\medskip
\roster
\item "{$(I)$}"  If $(p_1,r_1)\in A_{\mu^+}$, $\beta^*\in u^{r_1}_0$ then
there is $(p_2,r_2)\in A_{\mu^+}$ such that
$$
(p_1,r_1)\leq (p_2,r_2),\quad \eta^{p_2}_{\beta^*,0}\in v^{r_2}\quad\text{
and }\quad \eta^{p_2}_{\beta^*,1}=(\nu_{\eta^{p_2}_{\beta^*,0}})^{r_2}.
$$
\endroster
[Why? By clause $(A)$ above it is enough to find such $(p_2,r_2)$ in
$A^\prime_\beta$. Take $\beta^{**}\in u^{r_1}_0\cap\mu$ such that
$(\eta^{p_1}_{\beta^{**},0},\eta^{p_1}_{\beta^{**},1})=(\eta^{p_1}_{\beta^{*},
0},\eta^{p_1}_{\beta^{*},1})$, and let 
$$
\gamma=\sup\{\ell g(\nu)+1: \nu\in v^{r_1}\cup\{\eta^{p_1}_{\beta,\ell}:
\beta\in u^{p_1},\ell<2\}\}. 
$$ 
Define condition $(p_2,r_2)$ by: $\qquad u^{p_2}=u^{p_1}$, and
$$
\eta^{p_2}_{\beta,\ell} =\cases \eta^{p_1}_{\beta^*,0}\char 94 {\bar{\bold
0}}_\gamma\char 94\langle 1\rangle &\text{ if }
\qquad \beta \in \{\beta^*,\beta^{**}\},\ell=0, \\
  \eta^{p_1}_{\beta,1} \char 94 \langle 1\rangle &\text{ if } \qquad \beta \in 
\{ \beta^*,\beta^{**}\},\ell = 1, \\
  \eta^{p_1}_{\beta,\ell}\char 94 {\bar{\bold 0}}_{2\gamma+1} &\text{ if }
\qquad \beta \in u^{p_1}\setminus\{\beta^*,\beta^{**}\},\ell<2,
 \endcases
$$
$$
u^{r_2}_0=u^{r_1}_0,\quad v^{r_2}=v^{r_1}\cup\{\eta^{p_1}_{\beta^*,0}\char 94
{\bar{\bold 0}}_\gamma\char 94 \langle 1\rangle\},\quad\text{ and}
$$
$$
\nu^{r_2}_\rho=\cases\nu^{r_1}_\rho &\text{ if }\qquad\rho\in v^{r_1}_1,\\
   \eta^{p_1}_{\beta^*,1}\char 94\langle 1\rangle &\text{ if }\qquad
\rho=\eta^{p_1}_{\beta^*,0}\char 94{\bar{\bold 0}}_\gamma\char 94\langle
1\rangle.
 \endcases
$$
It is easy to check that $(p_2,r_2)\in A_{\mu^+}$ is as required.]
\medskip

Due to clause $(I)$ above we immediately get assertion $(a)$ of 1.3, i.e.~the
function $\eta\mapsto\bigcup\{\nu^r_\rho: \rho\vartriangleleft\eta,\ \rho\in
\bigcup\limits_{r\in G_{\underset\tilde {}\to R}} v^r\}$ is trivially
continuous and its range is contained in $^{\mu\geq} 2$, but clause $(I)$
implies that for each $\alpha\in{\underset\tilde {}\to {\Cal U}}$ the image of
${\underset\tilde {}\to \eta_{\alpha,0}}$ is ${\underset\tilde {}\to
\eta_{\alpha,1}\in {}^\mu 2}$ (and not a proper initial segment of it).  
\medskip

For clause $(b)$ assume that
\medskip
\roster
\item "{$(*)$}"  $\nu^* \in {}^{\mu >} 2$ and $\langle
\underset\tilde {}\to B_i:i < \mu \rangle$ is a sequence of $(P *
\underset\tilde {}\to R)$-names, \newline
$\Vdash_{P * \underset\tilde {}\to R} ``
\underset\tilde {}\to B_i \subseteq {}^\mu 2$ is nowhere dense and (for
simplicity) with no isolated points".
\endroster
\medskip

\noindent It suffices to prove
$$\Vdash_{P * \underset\tilde {}\to R}\text{``}\{{\underset\tilde {}\to
\eta_{\alpha,0}}:\alpha \in {\underset\tilde {}\to {\Cal U}} \text{ and }\nu^*
\triangleleft {\underset\tilde {}\to \eta_{\alpha,0}}\}\nsubseteq\bigcup_{i <
\mu}\underset\tilde {}\to B_i\text{''}.$$
Let $\underset\tilde {}\to T_i = \{\eta \restriction \gamma:\eta \in
\underset\tilde {}\to B_i$ and $\gamma < \mu\}$, so $\Vdash_{P *
\underset\tilde {}\to R}``\underset\tilde {}\to T_i \subseteq {}^{\mu >} 2$
is a nowhere dense subtree and $(\forall \eta \in {}^\mu 2)(\eta \in
\underset\tilde {}\to B_i \Rightarrow \dsize \bigwedge_{\gamma < \mu} \eta
\restriction \gamma \in {\underset\tilde {}\to T_i})$'', where
``$T \subseteq {}^{\mu >}2$ is a nowhere dense subtree" means:
$\langle \rangle \in T$, $(\forall \eta \in T)(\forall \gamma < \ell
g(\eta))(\eta \restriction \gamma \in T)$ and for every $\eta \in T$ and
$\gamma < \mu$ there are two $\vartriangleleft$-incomparable sequences
$\nu_0,\nu_1$ such that 
$$
\eta\triangleleft \nu_0 \in T \and \eta \triangleleft \nu_1 \in {}^{\mu >}2
\backslash T \and \ell g(\nu_0) \ge \gamma \and \ell g (\nu_1) \ge \gamma.
$$
We can find $\alpha < \mu^+$ such that $\alpha > \mu$ and for every $i < \mu$,
$\underset\tilde {}\to T_i$ is an $A_\alpha$-name (remember $\mu^{<\mu}=\mu$).
Let $(p_0,r_0) \in P *\underset\tilde {}\to R$. By the statement $(B)$ above
there is $(p_1,r_1)$ such that $(p_0,r_0) \le (p_1,r_1) \in A_{\mu^+}$.  We
can easily find $(p_2,r_2)$ such that $(p_1,r_1) \le (p_2,r_2) \in A_{\mu^+}$
and $\beta^* \in (\alpha,\mu^+)$ such that $\beta^* \in u^{r_2}_1$ and $\nu^*
\triangleleft \eta^{p_2}_{\beta^*,0}\in v^{r_2}_1$,
$(\nu_{\eta^{p_2}_{\beta^*,0}})^{r_2}=\eta^{p_2}_{\beta^*,1}$ (apply clauses
$(H)$ and $(I)$).

It suffices to prove that for each $i < \mu$, $(p_2,r_2) \Vdash
``{\underset\tilde {}\to \eta_{\beta^*,0}} \notin \underset\tilde {}\to B_i$". 
So let $i<\mu$ and $(p_2,r_2) \le (p_3,r_3) \in A_{\mu^+}$. Since $(p_3,r_3)
\restriction\alpha \in A_\alpha$, $\underset\tilde {}\to T_i$ is an
$A_\alpha$-name and $\Vdash_{A_\alpha} ``\underset\tilde {}\to T_i \subseteq
{}^{\mu >} 2$  is nowhere dense", there is a condition $(p_4,r_4) \in
A_\alpha$ such that $(p_3,r_3)\restriction \alpha \le (p_4,r_4)$ and for some
$\nu^{**}$ we have
$$
\eta^{p_3}_{\beta^*,0}\triangleleft\nu^{**} \in {}^{\mu >} 2\quad\text{ and }
\quad(p_4,r_4) \Vdash_{A_\alpha}\text{``}\nu^{**} \notin \underset\tilde {}\to
T_i\text{''}.
$$ 
Lastly we define $(p_5,r_5) \in A'_{\mu^+}$ such that $(p_5,r_5)$ is above
$(p_3,r_3)$ and above $(p_4,r_4)$ and
$(p_5,r_5)\Vdash$``$v^{**}\vartriangleleft\eta_{\beta^*,0}$''. So let
$$
u^{p_5} = u^{p_3} \cup u^{p_4},
$$

$$
\eta^{p_5}_{\beta,\ell} =\cases \eta^{p_4}_{\beta,\ell} &\text{ if }
\qquad \beta \in u^{p_4}, \\
  \eta^{p_3}_{\beta,\ell} &\text{ if } \qquad \beta \in u^{p_3} \setminus
(u^{p_4} \cup \{ \beta^* \}), \\ 
  \nu^{**} &\text{ if } \qquad \beta=\beta^*,\ell = 0, \\
  \bigcup\{\nu^{r_4}_\rho:\rho\in v^{r_4}\text{ and }\rho\trianglelefteq
  \nu^{**}\} &\text{ if } \qquad \beta=\beta^*,\ell = 1,  
\endcases
$$
$$
r_5=(u_0^{r_4}\cup u_0^{r_3},v^{r_4},\bar{\nu}^{r_4}).
$$
By clause $(A)$ there is $(p_6,r_6)$ such that $(p_5,r_5) \le (p_6,r_6) \in
A_{\mu^+}$. This clearly suffices. Thus we finish proving clause $(b)$ of the
conclusion of 1.3. 

Together we get all the required conclusions. \hfill$\square_{1.3}$
\enddemo
\bigskip

\proclaim{1.4 Theorem}  Assume that $\mu = \mu^{< \mu},2^\mu = \mu^+,
2^{\mu^+} = \mu^{+2}$ and $S$ is a stationary subset of  
$$
S^{\mu^{+2}}_{\mu^+} =: \{ \delta <\mu^{+2}:\text{cf}(\delta) = \mu^+\}.
$$
Suppose that $\diamondsuit_S$ holds. 

\underbar{Then} there is a $\mu$-closed, $\mu^+$-c.c. forcing notion $P$ of
cardinality $\mu^{+2}$ such that, in $V^P$, 
\medskip
\roster
\item "{$(*)$}"  for every function $f:{}^\mu 2 \rightarrow {}^\mu 2$, there
is a non $\mu$-meagre set $A \subseteq {}^\mu 2$ such that $f\restriction A$
is continuous. 
\endroster
\endproclaim
\bigskip

\demo{Proof}  Without loss of generality we may also assume that
$S^{\mu^{+2}}_{\mu^+}\backslash S$ is stationary.

For $\alpha \le \mu^{+2}$ let ${\frak K}_\alpha$ be the family of
sequences $\bar Q = \langle P_i,\underset\tilde {}\to Q_i:i < \alpha
\rangle$ such that
\medskip
\roster
\item "{$(a)$}"  $\bar Q$ is iteration with support $< \mu$ and $P_\alpha$ is 
the limit,
\item "{$(b)$}"  for $i\in\alpha\cap\mu^{+2}\backslash S$, ${\underset\tilde
{}\to Q_i}=Q_i$ is $({}^{\mu >}2, \triangleleft)$ (i.e. a $\mu$-Cohen forcing
notion) and we denote by ${\underset\tilde {}\to\eta^{\bar Q}_{i+1}}$ the name
of a Cohen generic element of ${}^\mu 2$ adjoined by $Q_i$,
\item "{$(c)$}"  if $\beta_1 < \beta_2 \le \alpha$, $\beta_2<\mu^{+2}$,
$\beta_1 \notin S$ then $P_{\beta_2}/P_{\beta_1}$ and $P_{\beta_2} /
(P_{\beta_1} * Q_{\beta_1})$ are $\mu$-Cohen forcing notions (and for clarity
we fix representation of $P_{\beta_2}$ as
$$
P_{\beta_1}\times Q_{\beta_1}\times\prod\{ Q^{\beta_1,\beta_2}_j:j <
j(\beta_1,\beta_2)\}
$$  
with support $< \mu$),
\item "{$(d)$}" we may assume that each ${\underset\tilde {}\to Q_i}$ (for
$i<\alpha$) is $\mu$-closed, $\mu^+$-cc and of cardinality $\leq\mu^+$. 
\endroster
\medskip

\noindent Note that by clause $(c)$, each $P_\beta$ is a $\mu$-closed forcing
notion satisfying the $\mu^+$-cc.

\noindent
Let ${\frak K} = \dsize \bigcup_{\alpha < \mu^{+2}} {\frak K}_\alpha$.

There is a natural ordering of ${\frak K}$:
$$
\bar Q^1\leq\bar Q^2\qquad\text{ if and only if }\qquad\bar Q^1 = \bar Q^2
\restriction \ell g (\bar Q^1).
$$
This partial order is $\mu^{+2}$-complete and every (strictly) increasing
chain of length $\mu^{+2}$ has a limit in ${\frak K}_{\mu^{+2}}$. 

Note that by clause $(c)$ we have that
\medskip
\roster
\item "{$(\boxtimes)$}"  if $\bar Q \in {\frak K}_\alpha,\beta<\alpha,\beta
\notin S$ and $\Vdash_{P_\beta} ``\underset\tilde {}\to B \subseteq {}^\mu 2$
is not $\mu$-meagre"

then $\Vdash_{P_\alpha} ``\underset\tilde {}\to B \subseteq {}^\mu 2$ is not
$\mu$-meagre".
\endroster
\medskip

\noindent
Because of this, and by $\diamondsuit_S$, it is enough to show that 
\medskip
\roster
\item "{$(\bigotimes)$}"  \underbar{if} $\bar Q\in{\frak K}_{\mu^{+2}}$, $P=
P_{\mu^{+2}} = \text{ lim}(\bar Q)$, and $\underset\tilde {}\to f$ is a
$P$-name of a function from ${}^\mu 2$ to ${}^\mu 2$ \underbar{then} for some
club $E$ of $\mu^{+2}$, for every $\delta \in S \cap E$ we have
{\roster
\itemitem{ $(\alpha)$ }  $i < \delta \Rightarrow \underset\tilde {}\to f
({\underset\tilde {}\to \eta^{\bar{Q}}_{i+1}})$ is a $P_\delta$-name,
\itemitem{ $(\beta)$ }  there is a $\bar Q^{\delta + 1} \in {\frak K}
_{\delta + 1}$ such that $\bar Q^{\delta + 1} \restriction \delta = \bar Q
\restriction \delta$ and 

$\Vdash_{\text{lim}(\bar Q^{\delta + 1})}$ ``for some $A \subseteq \delta,|A|
= \mu^+$, we have that

$\qquad \qquad (i) \quad
\{ {\underset\tilde {}\to \eta^{\bar Q}_{i+1}}:i \in A\}$ is not meagre,
\newline

$\qquad \qquad (ii) \quad$
the function $\{( {\underset\tilde {}\to \eta^{\bar Q}_{i+1}},
\underset\tilde {}\to f({\underset\tilde {}\to \eta^{\bar Q}_{i+1}})):i \in
A\}$ \newline

$\qquad \qquad \qquad \quad$ is continuous".
\endroster}
\endroster
(Note: ${\underset\tilde {}\to \eta^{\bar Q}_{i+1}}$ is a $P_\delta$-name 
(as $i+1 < \delta$) and
also $\underset\tilde {}\to f({\underset\tilde {}\to \eta^{\bar Q}_{i+1}})$
is a $P_\delta$-name (by clause $(\alpha)$)).
\enddemo
\bigskip

\demo{Proof of $(\bigotimes)$} For each $\delta \in S^{\mu^{+2}}_{\mu^+}
\backslash S$, $\underset\tilde {}\to f({\underset\tilde {}\to \eta^{\bar
Q}_{\delta+1}})$ is a $P_{\mu^{+2}}$-name of a member of ${}^\mu 2$, so for
some $\gamma_\delta>\delta + 1$, it is a $P_{\gamma_\delta}$-name and we can
demand that $\gamma_\delta$ is a successor ordinal. As 
$$
P_{\gamma_\delta}=P_\delta\times Q_\delta\times\prod_{j<j(\delta,
\gamma_\delta)} Q^{\delta,\gamma_\delta}_j
$$
(with support $<\mu$), we find a set $w_\delta\subseteq
j(\delta,\gamma_\delta)$ of cardinality $\le \mu$ such that $\underset\tilde
{}\to f({\underset\tilde {}\to \eta^{\bar Q}_{\delta+1}})$ is a $P_\delta\times
Q_\delta \times\prod\limits_{j \in w_\delta} Q^{\delta,\gamma_\delta}_j$-name.
Also as cf$(\delta) = \mu^+$ there is $\beta_\delta<\delta$ such that
$\underset\tilde {}\to f({\underset\tilde {}\to\eta^{\bar Q}_{\delta+1}})$ is
a $P_{\beta_\delta} \times Q_\delta \times \prod\limits_{j\in w_\delta}
Q^{\delta,\gamma_\delta}_j$-name and we can demand that $\beta_\delta$ is a
successor ordinal.  By Fodor's lemma (as $2^\mu = \mu^+$) for some stationary
$S_1 \subseteq S^{\mu^{+2}}_{\mu^+} \backslash S$ (and $\zeta < \mu^+$,
$\beta^* < \mu^{+2}$) we have that
\medskip
\roster
\item "{$(\alpha)$}"  $\delta \in S_1 \Rightarrow \beta_\delta = \beta^*$
\item "{$(\beta)$}"  $\delta \in S_1 \Rightarrow \text{ otp}(w_\delta) =
\zeta$
\item "{$(\gamma)$}"  for $\delta_1,\delta_2 \in S_1$ we demand
$f_{\delta_2,\delta_1}(\underset\tilde {}\to f({\underset\tilde {}\to
\eta^{\bar Q}_{\delta_1+1}})) = \underset\tilde {}\to f({\underset\tilde {}\to
\eta^{\bar Q}_{\delta_2+1}})$, where $f_{\delta_2,\delta_1}$ is the natural
isomorphism  
$$
\text{from }\quad P_{\beta^*} \times Q_{\delta_1} \times \prod_{j \in
w_{\delta_1}} Q^{\delta_1,\gamma_{\delta_1}}_j\quad\text{ onto }\quad 
P_{\beta^*} \times Q_{\delta_2} \times \prod_{j \in w_{\delta_2}} Q^{\delta_2,
\gamma_{\delta_2}}_j.
$$
\endroster
\medskip

\noindent
Clearly $Q^\delta=\prod\limits_{j\in w_\delta} Q^{\delta,\gamma_\delta}_j$ is
isomorphic to Cohen forcing.

Lastly, let
$$
E = \biggl\{ \delta < \mu^{+2}:\delta = \sup(S_1 \cap \delta) \text{ and }
  \delta_1 \in S_1 \cap \delta \Rightarrow \gamma_{\delta_1} < \delta
\biggr\}.
$$
Clearly the set $E$ is a club of $\mu^{+2}$. We are going to show that for
each $\delta^* \in E \cap S$ there is $\bar Q^{\delta^*+1}$ as required in
$(\bigotimes)$.  

Let $\delta^*\in E\cap S$. We can find an increasing continuous sequence
$\langle\beta_\varepsilon:\varepsilon<\mu^+\rangle$  with limit $\delta^*$ and
a sequence $\langle \delta_\varepsilon:\varepsilon < \mu^+ \rangle$ such that
$\beta_0 = \beta^*$, $\beta_\varepsilon$ a successor ordinal if $\varepsilon$
is a non-limit, $\beta_\varepsilon\notin S$, $\delta_\varepsilon\in S_1$, and
$\beta_\varepsilon<\delta_\varepsilon<\gamma_{\delta_\varepsilon}<
\beta_{\varepsilon+1}$. Let ${\underset\tilde {}\to\eta^* _{\varepsilon,0}}$
and ${\underset\tilde {}\to \eta^*_{\varepsilon,1}}$ be
$P_\delta/P_{\beta_0}$-names such that ${\underset\tilde {}\to \eta^*
_{\varepsilon,0}}$ is ${\underset\tilde {}\to\eta^{\bar Q}_{\delta_\varepsilon
+1}}$ and ${\underset\tilde {}\to\eta^*_{\varepsilon,1}}$ is a name for the
Cohen subset of $\mu$ added by $Q^{\delta_\varepsilon}$.  

Now let
$$
R^0_{\delta^*}=:\prod[\{Q_{\delta_\varepsilon}:\varepsilon < \mu^+\} \cup \{ 
Q^{\delta_\varepsilon}:\varepsilon < \mu^+\}]
$$
and note that $R^0_{\delta^*}\lesdot P_{\delta^*}/P_{\beta^*}$. The forcing
notion $R^0_{\delta^*}$ is naturally isomorphic to $P$ from 1.3 with
${\underset\tilde {}\to \eta^*_{\varepsilon,\ell}}$ corresponding to 
${\underset\tilde {}\to\eta_{\varepsilon,\ell}}$. Moreover the quotient
$P_{\delta^*}/(P_{\beta^*}\times R^0_{\delta^*})$ is a $\mu$-Cohen forcing
notion. [Why? For each $\varepsilon<\mu^+$, $\delta_\varepsilon,
\gamma_{\delta_\varepsilon},\beta_\varepsilon\notin S$ and hence we may write
$P_{\beta_{\varepsilon+1}}$ as 
$$
P_{\beta_\varepsilon}\times Q_{\beta_\varepsilon}
\times\prod_{j<j(\beta_\varepsilon,\delta_\varepsilon)} Q^{\beta_\varepsilon,
\delta_\varepsilon}_j \times Q_{\delta_\varepsilon}\times
\prod_{j<j(\delta_\varepsilon,\gamma_{\delta_\varepsilon})}
Q^{\delta_\varepsilon,\gamma_{\delta_\varepsilon}}_j\times
Q_{\gamma_{\delta_\varepsilon}} \times \prod_{j<j(\gamma_{\delta_\varepsilon},
\beta_{\varepsilon+1})}Q^{\gamma_{\delta_\varepsilon},
\beta_{\varepsilon+1}}_j.
$$
But
$$
\prod_{j<j(\delta_\varepsilon,\gamma_{\delta_\varepsilon})}
Q^{\delta_\varepsilon,\gamma_{\delta_\varepsilon}}_j = Q^{\delta_\varepsilon}
\times \prod_{j\in j(\delta_\varepsilon,\gamma_{\delta_\varepsilon})\setminus
w_{\delta_\varepsilon}} Q^{\delta_\varepsilon,\gamma_{\delta_\varepsilon}}_j.
$$
So we may represent $P_{\beta_{\varepsilon+1}}$ as
$$
P_{\beta_{\varepsilon+1}}=P_{\beta_\varepsilon}\times Q_{\delta_\varepsilon}
\times Q^{\delta_\varepsilon}\times R^*_\varepsilon,
$$
where $R^*_\varepsilon$ is a product (with support $<\mu$) of $\mu$-Cohen
forcing notions. Now, since the sequence $\langle\beta_\varepsilon:\varepsilon
<\mu^+\rangle$ is increasing and continuous with limit $\delta^*$ we may write
$$
P_{\delta^*}=P_{\beta^*}\times\prod_{\varepsilon<\mu^+} Q_{\delta_\varepsilon}
\times\prod_{\varepsilon<\mu^+} Q^{\delta_\varepsilon}
\times\prod_{\varepsilon<\mu^+} R^*_\varepsilon
$$
(all products with $<\mu$ support).]

Let ${\underset\tilde {}\to Q_{\delta^*}}$ be $R[\langle {\underset\tilde
{}\to \eta^*_{\varepsilon,0}}, {\underset\tilde {}\to\eta^*_{\varepsilon,1}}:
\varepsilon < \mu^+ \rangle]$ (from 1.3), $\bar Q^{\delta^*+1} = (\bar Q
\restriction \delta^*) \char 94 \langle (P_{\delta^*},{\underset\tilde {}\to
Q_{\delta^*}}) \rangle$. To check that $\bar Q^{\delta^*+1}$ is as desired we
use 1.3. Thus $P_{\delta^*}/P_\beta$ (for $\beta<\delta$, $\beta\notin S$)
is $\mu$-Cohen by 1.3(c) and because, as said above, $P_{\delta^*}/
(P_{\beta^*}\times R^0_{\delta^*})$ is $\mu$-Cohen. It follows from 1.3(a)
that the function ${\underset\tilde {}\to \eta^*_{\varepsilon,0}}\mapsto
{\underset\tilde {}\to\eta^*_{\varepsilon,1}}$ (for $\varepsilon\in
{\underset\tilde{}\to{\Cal U}}$) is continuous and hence, as the names
${\underset\tilde {}\to f}({\underset\tilde {}\to\eta^{\bar
Q}_{\delta_\varepsilon+1}})$ are isomorphic, the function ${\underset\tilde
{}\to \eta^{\bar Q}_{\delta_\varepsilon+1}}\mapsto {\underset\tilde{}\to
f}({\underset\tilde {}\to\eta^{\bar Q}_{\delta_\varepsilon+1}})$ (for
$\varepsilon \in{\underset \tilde{}\to {\Cal U}}$) is continuous (it is
obviously Borel, which actually suffices, but in fact, as the ${\underset
\tilde{}\to\eta^*_{\varepsilon,0}}$ are Cohen over the definition of
${\underset \tilde{}\to f}$ it is continuous). Finally, by
1.3(b), 
$$
\Vdash_{P_{\delta+1}}\text{``}\{{\underset\tilde{}\to\eta^{\bar
Q}_{\delta_\varepsilon+1}}: \varepsilon\in{\underset\tilde{}\to{\Cal U}}\}
\text{ is not meagre''.}
$$
This finishes the proof. \hfill$\square_{1.4}$
\enddemo
\bigskip

\remark{1.5 Remark:} 1. Alternatively, one may force the iteration instead of
using $\diamondsuit_S$.

\noindent 2. The forcing notion from 1.4 is constructed like ``$\lambda$-free
not free (Boolean) algebra of cardinality $\lambda$'', with {\it free}
interpreted as {\it $\mu$-Cohen}.
\endremark

\head {\S2 Final remarks} \endhead

\noindent
Haim Judah asked us whether adding Cohen reals is enough.
\bigskip

\proclaim{2.1 Observation}  If  $P$  is the forcing notion for adding
$2^{\aleph_0}$ Cohen real, then in  $V^P$ there is $f:{\Bbb R} \rightarrow
{\Bbb R}$ which is not continuous on any uncountable set.
\endproclaim
\bigskip

\demo{Proof}  Let $\lambda = 2^{\aleph_0}$, $P = \{ f:f$  is a finite function
from $\lambda$ to $\{0,1\}\}$. For any set $w\subseteq\lambda$ let $P_w=\{f\in
P:\text{Dom}(f) \subseteq w\}$.  It is known that $P_w\lesdot P = P_\lambda$
and for any $P$-name ${\underset\tilde {}\to \eta}$ of a real (i.e. $\Vdash_P$
``${\underset\tilde {}\to \eta} \in {}^\omega 2$'') for some countable
$w\subseteq \lambda$, ${\underset\tilde {}\to \eta}$ is a $P_w$-name.  For
every ordinal $\alpha < \lambda$ let ${\underset\tilde {}\to r_\alpha} \in
{}^\omega 2$ be defined by: 
$$
{\underset\tilde {}\to r_\alpha}(n) = \ell\text{ if and only if\ \ \ for some
}f \in {\underset\tilde {}\to G_P},\ f(\omega \alpha+n)=\ell.
$$  
Clearly ${\underset\tilde {}\to r_\alpha}$ is a Cohen real. Let $\langle
{\underset\tilde {}\to s_\zeta}:\zeta<\lambda\rangle$ list all $P$-names of
reals, so for some countable  $w_\zeta \subseteq \lambda$,
${\underset\tilde {}\to s_\zeta}$ is a $P_{w_\zeta}$-name. Without loss of
generality  
$$
\alpha + \omega = \beta + \omega \and \alpha \in  w_\zeta\qquad \Rightarrow
\qquad \beta\in w_\zeta.
$$  
We now define by induction on $\zeta<\lambda$ an ordinal $\alpha_\zeta <
\lambda$  as follows:
$$
\alpha_\zeta \text{ is the first ordinal }  \alpha \text{ such that }
\alpha\notin\{\alpha_\varepsilon:\varepsilon < \zeta\}\cup\{\beta:\omega+k
\beta\in w_\zeta\text{ for some }k\in\omega\}.
$$
Let $G \subseteq  P$ be a generic filter over $V$. In $V[G]$ we define a
function $F:{}^\omega 2 \rightarrow {}^\omega 2$ as follows:
$$
\text{for }s \in {}^\omega 2\ \ \text{ let }\ \gamma(s) = \min\{\zeta:s =
{\underset\tilde {}\to s_\zeta}[G]\}\ \text{ and }\ F(s) = {\underset\tilde
{}\to r_{\alpha_{\gamma(s)}}}.
$$
Suppose now that, in $V[G]$, a set $A \subseteq {}^\omega 2$  is uncountable
and $F\restriction A$  is continuous. Thus we have a Borel function
$F':{}^\omega 2\longrightarrow {}^\omega 2$ such that $F \restriction A = F'
\restriction A$. The function $F'$ is coded by a single real $s \in {}^\omega
2$, $\gamma(s)=\zeta$. 

Let $p_0 \in G$, ${\underset\tilde {}\to A},{\underset\tilde {}\to F},
{\underset\tilde {}\to F'}, {\underset\tilde{}\to s}$ and $\zeta<\lambda$ be
such that 
$$
p_0\Vdash_P\text{``}{\underset\tilde {}\to A},{\underset\tilde {}\to F},
{\underset\tilde {}\to F'},{\underset\tilde{}\to s}={\underset\tilde{}\to
s_\zeta},\zeta=\gamma({\underset\tilde{}\to s})\text{ are as above''}.
$$
The set $\{\xi<\lambda: [\omega\alpha_\xi,\omega\alpha_\xi+\omega)\cap w_\zeta
\neq\emptyset\}$ is countable, $p \Vdash_P$``${\underset\tilde {}\to A}$ is
uncountable'', so we find $\xi<\lambda$ and $p_1\geq p_0$, $p_1\in G$ such
that
$$
p_1\Vdash_P\text{``}{\underset\tilde {}\to s_\xi}\in {\underset\tilde {}\to A}
\text{ and } \gamma({\underset\tilde{}\to s_\xi})=\xi\text{''}
$$
and $[\omega\alpha_\xi,\omega\alpha_\xi+\omega)\cap w_\zeta=\emptyset$. Now
$$
p_1\Vdash_P\text{``}{\underset\tilde{}\to F'}({\underset\tilde {}\to s_\xi}) =
{\underset\tilde{}\to F}({\underset\tilde {}\to s_\xi}) =
{\underset\tilde{}\to r_{\alpha_\xi}}\text{''}.
$$
But we can compute ${\underset\tilde{}\to F'}({\underset\tilde {}\to s_\xi})$
in $V[G\cap P_{w_\zeta\cup w_\xi}]$ and (by the choice of $\xi$)
${\underset\tilde{}\to r_{\alpha_\xi}}$ is Cohen over  $V[G\cap P_{w_\zeta\cup
w_\xi}]$ (remember the choice of $\alpha_\xi$'s), which is a contradiction.
\hfill$\square_{2.1}$ 
\enddemo
\bigskip

\remark{2.2 Remark}  Instead of forcing it suffices to assume the existence of
a Luzin set in a strong sense. 
\endremark
\bigskip

\centerline {$* \qquad * \qquad *$}
\bigskip

We call $A \subseteq {}^\omega 2$ a \underbar{Luzin} set if it is uncountable
and non-meagre and any uncountable $B \subseteq A$ is non-meagre.

\definition{2.3 Definition}  A forcing notion $Q$ is Luzin-preserving
\underbar{if} it satisfies the c.c.c. and
\medskip
\roster
\item "{$(*)$}"  for any  Luzin set $A\subseteq {}^\omega 2$ and a sequence
$\langle p_\eta:\eta \in  A \rangle$  of members of  $Q$, we have that

\item "{$(**)$}"  for some  $\eta  \in  A$ \newline
$p_\eta \Vdash_Q$ ``for some $\rho \in {}^{\omega >} 2$, \newline
$(\forall \nu)[\rho \trianglelefteq \nu \in {}^{\omega >} 2 \Rightarrow$
for uncountably many $\eta_1 \in A,\nu \triangleleft \eta_1 \and p_\eta \in  
{\underset\tilde {}\to G_Q}$]".
\endroster
\enddefinition
\bigskip

\proclaim{2.4 Claim}  0) Assume there are Luzin sets, then in Definition
2.3 ``it satisfies the c.c.c." follows from $(*)$ of 2.3. \newline
1)  A forcing notion  $Q$  is Luzin-preserving \underbar{iff} it
satisfies the c.c.c. 
and for any Luzin  $A \subseteq {}^\omega 2$  we have, $\Vdash_Q$ ``$A$  is 
Luzin". \newline
2)  Being Luzin-preserving is preserved by composition. \newline
3)  If  $\langle P_i:i < \delta \rangle$  is an increasing continuous sequence
of Luzin-preserving forcing notions \underbar{then}
$\dsize \bigcup_{i < \delta} P_i$ is a Luzin-preserving forcing notions.
\newline
4)  The statement ``a forcing notion $P$ is Luzin preserving" depend on
$P$ only up to equivalence of forcing notions \newline
5)  In Definition 2.3 we can weaken the demand to:
\medskip
\roster
\item "{$(*)'$}"  for some $\eta \in A$ there is $q$ such that
$p_\eta \le q \in Q$ and $q \Vdash_Q$ ``...".
\endroster
\endproclaim
\bigskip

\demo{Proof}  0)  Should be clear. \newline
1) $\Rightarrow$ \underbar{(the ``only if" implication)}. \newline By
definition, $Q$ satisfies the c.c.c. Assume $A$ is Luzin but $p \Vdash ``A$ is
not Luzin". Then for some $Q$-names ${\underset\tilde{}\to\eta_i}$ (for
$i<\omega_1)$, we have 
$$
p \Vdash\text{``}{\underset\tilde {}\to \eta_i}\in A\text{ for }i<\omega_1,\
\bigwedge_{i < j} {\underset\tilde {}\to \eta_i}\neq {\underset\tilde {}\to 
\eta_j}\text{ and }\{ {\underset\tilde {}\to\eta_i}: i<\omega_1\}\text{ is not
meagre''}.
$$
Thus changing the ${\underset\tilde {}\to \eta_i}$'s we can get
\medskip
\roster
\item "{$(\bigoplus)$}"  $p \Vdash ``{\underset\tilde {}\to \eta_i} \in A$  for
$i < \omega_1$,  $\dsize \bigwedge_{i<j} {\underset\tilde {}\to \eta_i} \ne
{\underset\tilde {}\to \eta_j}$ and  $\{ {\underset\tilde {}\to \eta_i}:
i < \omega_1 \}$  is nowhere dense".
\endroster
Let $p \le p_i \in Q$, $p_i \Vdash ``{\underset\tilde {}\to \eta_i} = \nu_i"$
so $\nu_i \in \big({}^\omega 2 \big)^V$ and necessarily $\nu_i \in A$.  For no
$\nu$ the set $\{i:\nu_i = \nu \}$ is uncountable as then $Q$ fails the c.c.c.
So without loss of generality $i\ne j \Rightarrow \nu_i \ne \nu_j$.  Let $A' =
\{ \nu_i:i < \omega_1 \}$. It is an uncountable subset of $A$ and hence it is
Luzin. Let $p_\nu = p_i$ if $\nu = \nu_i$.  Apply $(*)$ of Definition 2.3 to
$\langle p_\nu:\nu \in A' \rangle$ and get contradiction to $(\bigoplus)$
above.
\bigskip

\noindent
$\Leftarrow$ (\underbar{the ``if'' implication}). Of course $Q$ satisfies the
c.c.c.  Let $\langle p_\eta:\eta \in A \rangle$ be given, $A \subseteq
{}^\omega 2$ a Luzin set and $p_\eta \in Q$ (for $\eta\in A$) and we should
show that $(**)$ of Definition 2.3 holds.  Then ${\underset\tilde {}\to A}=\{
\eta \in A:p_\eta \in {\underset\tilde{}\to G_Q}\}$ is a $Q$-name of a
subset of $A$.  It is known that for all but countably many $\eta \in A$,
we have 
$$
p_\eta\Vdash``{\underset\tilde {}\to A} \subseteq  A \text{ is uncountable''}
$$
so that (by the assumption of the implication we are proving)
$$
p_\eta \Vdash ``{\underset\tilde {}\to A} \text{ is non-meagre, hence
nowhere meagre above some } {\underset\tilde {}\to \nu_\eta} \in
{}^{\omega >} 2".
$$
\medskip

\noindent
We have really finished, but we can elaborate.  As we can replace $A$ by $A
\backslash B$ for any countable $B \subseteq A$. Without loss of generality
this holds for every $\eta \in A$.  Now for $\eta \in A$ and
$p'_\eta,p_\eta \le p'_\eta \in P$ there are $\nu^\ast_\eta$ and $q_\eta$,
$p'_\eta \le q_\eta \in Q$ such that $q_\eta
\Vdash_Q ``{\underset\tilde {}\to \nu_\eta} = \nu^\ast_\eta$".  So for all
$\eta \in A$ we have 
$$
q_\eta\Vdash\text{``}\{\nu\in A:p_\nu\in{\underset\tilde{}\to G_Q}\}\text{ is
uncountable and everywhere non-meagre above }\nu^*_\eta\text{''}.
$$
As this holds for any $p'_\eta$, $p_\eta \le p'_\eta \in P$ we are done.
\newline 
2) Follow by 2.4(1) (and as c.c.c. is preserved). \newline 3) Like the
preservation of c.c.c. \newline 4) Easy (e.g. use 2.4(1)). \newline 5) Left to
the reader. \hfill$\square_{2.4}$
\enddemo
\newpage

\newpage
\bigskip
REFERENCES
\bigskip

\bibliographystyle{lit-plain}
\bibliography{lista,listb,listx}

\shlhetal

\enddocument
\bye